\documentclass[12pt]{article}
\usepackage{amsmath,amsfonts,latexsym,amsthm}
\newcommand{\F}{\mathbf F_q}
\newcommand{\K}{\overline{K}_c}
\begin{document}
\newtheorem{lem}{Lemma}
\newtheorem{teo}{Theorem}
\newtheorem*{cor}{Corollary}
\pagestyle{plain}
\title{$\F$-Linear Calculus over Function Fields}
\author{Anatoly N. Kochubei\\ \footnotesize Institute of Mathematics,\\
\footnotesize Ukrainian National Academy of Sciences,\\
\footnotesize Tereshchenkivska 3, Kiev, 252601 Ukraine
\\ \footnotesize E-mail: \ ank@ank.kiev.ua}
\date{}
\maketitle
\vspace*{2cm}
Copyright Notice: This work has been submitted to Academic Press for
possible publication. Copyright may be transferred without notice,
after which this version may no longer be accessible.
\newpage
\vspace*{8cm}
\begin{abstract}
We define analogues of higher derivatives for $\F$-linear functions over
the field of formal Laurent series with coefficients in $\F$. This
results in a formula for Taylor coefficients of a $\F$-linear
holomorphic function, a definition of classes of $\F$-linear smooth
functions which are characterized in terms of coefficients of
their Fourier-Carlitz expansions. A Volkenborn-type integration
theory for $\F$-linear functions is developed; in particular, an
integral representation of the Carlitz logarithm is obtained.
\end{abstract}
\vspace{2cm}
\newpage
\section{INTRODUCTION}
It is well known that any non-discrete, locally compact field of
characteristic $p$ is isomorphic to the field $K$ of formal
Laurent series with coefficients from the Galois field $\F$,
$q=p^\gamma $, $\gamma \in \mathbf Z_+$. The foundations of
analysis over $K$ were laid in a series of papers by Carlitz
(see, in particular, [3-5]), Wagner [20, 21], Goss [6-8], and
Thakur [17, 18].

An interesting property of many functions introduced in the above
works as analogues of classical elementary and special functions,
is their $\F$-linearity. Recall that a function $f:\ K\to \K$ (where
$\K$ is a completion of an algebraic closure of $K$) is called $\F$-linear
if $f(t_1+t_2)=f(t_1)+f(t_2)$ and $f(\alpha t)=\alpha f(t)$ for any
$t,t_1,t_2\in K$, $\alpha \in \F$. In a similar way one can define a
notion of a  $\F$-linear function on any $\F$-subspace of $K$,
and also a notion of a $\F$-linear operator on a vector space
over $K$.

Tools of classical calculus are not sufficient to study behavior
of $\F$-linear functions. For example, if such a function $f$ is
differentiable then $f'(t)\equiv \mbox{const}$, and all the
higher derivatives vanish irrespective of possible properties of
$f$. In particular, one cannot reconstruct the Taylor
coefficients of a holomorphic function -- the classical formula
contains the expression $\frac{f^{(n)}(t)}{n!}$ where both the
numerator and denominator vanish.

Thanks to Carlitz, we know the correct analogue of the factorial.
A counterpart of a higher derivative is given in this paper. In
fact, its main ingredient, the difference operator $\Delta
^{(n)}$,
$$
\left( \Delta ^{(n)}f\right) (t)=\Delta ^{(n-1)}f(xt)-
x^{q^{n-1}}\Delta ^{(n-1)}f(t),\quad n\ge 1;\ \ \Delta ^{(0)}=
\mbox{id},
$$
where $x$ is a prime element in $K$, was introduced by Carlitz
and often used subsequently.

Our approach is based on identities for the difference operator
$\Delta =\Delta ^{(1)}$, which emerge as a function field
analogue of the Bargmann-Fock representation of the canonical
commutation relations of quantum mechanics [9, 14]. Note that an
analogue of the Schr\"odinger representation was obtained in
[12].

We will also use $\Delta ^{(n)}$ in order to characterize
``smoothness'' of $\F$-linear functions and obtain an exact
correspondence between the degree of smoothness (or analyticity)
and the decay rate for coefficients of Fourier-Carlitz
expansions. Note that in the $p$-adic case (where smoothness is
understood in a conventional sense) such a correspondence was
found in [1, 2, 15].

Having a new ``derivative'', we can introduce a kind of an
antiderivative, and a Volkenborn type integral (see [15, 19] for the
case of zero characteristic), which leads, in particular, to an
integral representation of the Carlitz logarithm, establishing a direct
connection between the latter and the Carlitz module operation.

\section{PRELIMINARIES}
\subsection{{\it Carlitz Basis} {\rm [3, 4, 6, 12, 20]}}

Denote by $|\cdot |$ the non-Archimedean absolute value on $K$; if
$z\in K$,
$$
z=\sum _{i=n}^\infty \zeta_ix^i,\quad n\in \mathbf Z,\ \zeta_i\in \F \ ,
\ \zeta_n\ne 0,
$$
then $|z|=q^{-n}$. It is well known that this valuation can be
extended onto $\K$. Let $O=\{ z\in K:\ |z|\le 1\}$ be the ring of
integers in $K$. The ring $\F [x]$ of polynomials (in the
indeterminate $x$) with coefficients from $\F$ is dense in $O$.

Let $C_0(O,\K )$ be the Banach space of all $\F$-linear continuous
functions $O\to \K$, with the supremum norm $\| \cdot \|$. Any
function $\varphi \in C_0(O,\K )$ admits a unique representation
as a uniformly convergent series
$$
\varphi =\sum \limits _{i=0}^\infty c_if_i,\quad c_i\in \K,\
c_i\to 0,
$$
satisfying the orthonormality condition
$$
\| \varphi \| =\sup_{i\ge 0}|c_i|.
$$
Here $\{ f_i\}$ is the sequence of the normalized Carlitz
$\F$-linear polynomials defined as follows.

Let $e_0(t)=t$,
\begin{equation}
e_i(t)=\prod \limits _{{m\in \mathbf F_q[x]\atop \deg m<i}}(t-m),\quad
i\ge 1.
\end{equation}
It is known [3, 6] that
$$
e_i(t)=\sum \limits _{j=0}^i(-1)^{i-j}\left[ i\atop j\right]
t^{q^j}\eqno (1')
$$
where
$$
\left[ i\atop j\right]=\frac{D_i}{D_jL_{i-j}^{q^j}},
$$
the elements $D_i,\ L_i\in K$ are defined as
$$
D_i=[i][i-1]^q\ldots [1]^{q^{i-1}};\ L_i=[i][i-1]\ldots [1]\
(i\ge 1);\ D_0=L_0=1,
$$
and $[i]=x^{q^i}-x\in O$. Finally,
$$
f_i(t)=D_i^{-1}e_i(t),\quad i=0,1,2,\ldots  .
$$

The basis $\{ f_i\}$ can be complemented up to an orthonormal
basis $\{ h_j\}$ of the space $C(O,\K )$ of all continuous
$\K$-valued functions on $O$, as follows. Let us write any
natural number $j$ as
\begin{equation}
j=\sum \limits _{i=0}^\nu \alpha _iq^i,\quad 0\le  \alpha _i<q,
\end{equation}
and set
$$
h_j(t)=\frac{G_j(t)}{\Gamma _j},\quad G_j(t)=\prod \limits
_{i=0}^\nu (e_i(t))^{\alpha _i},\quad  \Gamma _j=\prod \limits
_{i=0}^\nu D_i^{\alpha _i}.
$$
It is clear that $h_j$ is a polynomial of degree $j$, $h_{q^i}=f_i$.

\begin{lem}
The elements $D_m$, $L_m$, and $\Gamma _{q^m-1}$ are connected by
the identity
$$
\Gamma _{q^m-1}L_m=D_m.
$$
\end{lem}

{\it Proof}. By the definition,
$$
\Gamma _{q^m-1}=(D_0\ldots D_{m-1})^{q-1}.
$$
We have $D_m=\prod \limits_{i=1}^m[i]^{q^{m-i}}$,
$L_m=\prod \limits_{i=1}^m[i]$,
$$
D_0\ldots D_{m-1}=\prod \limits_{j=1}^{m-1}\prod
\limits_{i=1}^j[i]^{q^{j-i}}=\prod \limits_{i=1}^{m-1}\prod
\limits_{j=i}^{m-1}[i]^{q^{j-i}}
=\prod \limits_{i=1}^{m-1}[i]^{\sum \limits _{j=i}^{m-1}q^{j-i}}
=\prod \limits_{i=1}^{m-1}[i]^{\frac{q^{m-i}-1}{q-1}},
$$
so that
$$
\Gamma _{q^m-1}L_m=[m]\prod \limits_{i=1}^{m-1}[i]^{q^{m-i}}=D_m.
\qquad \Box
$$

Another system of polynomials introduced by Carlitz is given by
$$
g_j(t)=\prod \limits _{i=0}^\nu g_{\alpha _iq^i}(t),\quad
g_0(t)\equiv 1,
$$
where $j$ is expanded as in (2),
$$
g_{\alpha _iq^i}(t)=\left\{
\begin{array}{rl}
e_i^{\alpha _i}(t), & \mbox{if }\alpha _i<q-1,\\
e_i^{q-1}(t)-D_i^{q-1}, & \mbox{if }\alpha _i=q-1.
\end{array}\right.
$$

\begin{lem}
The system of polynomials
$$
\tau _m(t)=\frac{g_{q^m-1}(t)}{\Gamma _{q^m-1}},\quad m=0,1,2,\ldots
,\ t\in O,
$$
is orthonormal.
\end{lem}

{\it Proof}. By the definition,
$$
g_{q^m-1}(t)=\prod \limits _{i=0}^{m-1}\left( e_i^{q-1}(t)-
D_i^{q-1}\right) ,
$$
whence
$$
\tau _m(t)=\prod \limits _{i=0}^{m-1}\left( f_i^{q-1}(t)-1\right) ,
$$
so that $\| \tau _m\| \le 1$. Since $\tau _m(0)=(-1)^{m-1}$, we
have $|\tau _m(0)|=1$ and $\| \tau _m\| =1$.

In order to prove orthogonality, it suffices to show that for any
natural $m$, and any \linebreak $\lambda _1,\ldots ,\lambda _m\in \K$
\begin{equation}
\left\| \sum \limits _{k=0}^m\lambda _k\tau _k\right\| \ge
|\lambda _m|
\end{equation}
(see Proposition 50.4 in [15]).

It is known [6] that
$$
g_{q^m-1}(t)=\sum \limits _{j=0}^{q^m-1}G_j(t)g_{q^m-1-j}(0).
$$
It follows from the definition of $\Gamma _j$ that
$$
\Gamma _{q^m-1}=\Gamma _{q^m-1-j}\Gamma _j,\quad 0\le j\le q^m-1,
$$
so that
$$
\tau _m(t)=\sum \limits _{j=0}^{q^m-1}\sigma _{m,j}h_j(t),\quad
\sigma _{m,j}=\frac{g_{q^m-1-j}(0)}{\Gamma _{q^m-1-j}}.
$$

Now
$$
\sum \limits _{k=0}^m\lambda _k\tau _k(t)=\sum \limits _{j=0}^{q^m-1}
h_j(t)\sum \limits _{\log _q(j+1)\le k\le m}\lambda _k\sigma
_{k,j},
$$
$$
\left\| \sum \limits _{k=0}^m\lambda _k\tau _k\right\| \ge
|\lambda _m\sigma _{m,q^m-1}|,
$$
due to the orthonormality of $\{ h_j\}$. Since $\sigma _{m,q^m-
1}=\frac{g_0(0)}{\Gamma _0}=1$, we come to (3). $\quad \Box$

\subsection{{\it Canonical Commutation Relations} {\rm [12]}}

In the quantum mechanics of harmonic oscillator (see e.g. [16]) a
creation operator transforms a stationary state into a
stationary state of the next (higher) energy level, an annihilation
operator acts in the opposite way. In quantum field theory these
properties are used to obtain operators which change the number of
particles.

The analogues of the creation and annihilation operators (in the
Schr\"odinger representation) for our situation are as follows.
Let $R_q$ be a $\F$-linear operator on $C_0(O,\K )$ of the form
$R_qu=u^q$. Consider the operators
$$
a^+=R_q-I,\quad a^-=\sqrt[q]{\ }\circ \Delta
$$
where $I$ is the identity operator, $\Delta =\Delta ^{(1)}$ (see
Introduction), while $\sqrt[q]{\ }$ is the inverse to $R_q$. Both
operators $a^+$ and $a^-$ are $\F$-linear and continuous on
$C_0(O,\K )$ and
\begin{equation}
a^-a^+-a^+a^-=[1]^{1/q}I.
\end{equation}
The operator $a^+a^-$ possesses the orthonormal eigenbasis
$\{ f_i\} $,
\begin{equation}
(a^+a^-)f_i=[i]f_i,\quad i=0,1,2,\ldots ;
\end{equation}
$a^+$ and $a^-$ act upon the basis as follows:
\begin{equation}
a^+f_{i-1}=[i]f_i,\ \ a^-f_i=f_{i-1},\ i\ge 1;\ a^-f_0=0.
\end{equation}

We look at (4)-(6) as an analogue of the Schr\"odinger representation
(for an exposition see [14, 16] or other standard texts on
mathematical foundations of quantum mechanics) since in
non-Archimedean analysis spaces of continuous functions (or
sequence spaces isomorphic to them) are often used instead of
the Hilbert spaces of square integrable functions where the
operators of conventional quantum mechanics are defined. Note
also that a $p$-adic analogue of the Schr\"odinger representation
was found in [11].

It is interesting that the above operators possess remarkable
algebraic properties. The operator $\Delta $ is a derivation on
the ring of $\F$-linear functions from $O$ to itself, with
composition being the multiplication in the ring. The operator
$R_q-I$ is an $R_q$-derivation in the sense of [10].

The Bargmann-Fock representation is a realization of a structure
like (4)-(6) by operators on a space of holomorphic functions. A
simple construction for the case of a function field is given
below.

\section{HOLOMORPHIC FUNCTIONS AND THE \\ BARGMANN-FOCK
REPRESENTATION}

Consider a Banach space $H$ over the field $\K$ consisting of
$\F$-linear power series
\begin{equation}
u(t)=\sum \limits _{n=0}^\infty a_n\frac{t^{q^n}}{D_n},\quad
a_n\in \K,\ |a_n|\to 0.
\end{equation}
The norm in $H$ is given by
$$
\| u\| _A=\sup \limits _n|a_n|.
$$

Since
$$
|D_n|=q^{-1}\left( q^{-1}\right) ^q\ldots \left( q^{-1}\right)
^{q^{n-1}}=q^{-\frac{q^n-1}{q-1}}=\left( q^{-1/(q-1)}\right) ^{q^n-1},
$$
a series (7) defines a holomorphic function for $|t|\le q^{-1/(q-1)}$.

It is obvious that the sequence of the functions $\widetilde{f}_n(t)=
\frac{t^{q^n}}{D_n}$, $n=0,1,2,\ldots $, is an orthonormal basis of $H$.

The desired representation is given by the following operators on the
space $H$:
$$
\widetilde{a}^+=R_q,\quad \widetilde{a}^-=\sqrt[q]{\ }\circ \Delta
$$
(note that the form of the operators $\widetilde{a}^\pm$ is only
slightly different from the one of $a^\pm$, but they act on a
different Banach space!).

By a sraightforward computation based on the identities
$$
\Delta \left( t^{q^n}\right) =[n]t^{q^n},\quad
D_{n+1}=[n+1]D_n^q,
$$
we show that the relations (4)-(6) hold for the operators
$\widetilde{a}^\pm$, with $\widetilde{f}_n$ substituted for
$f_n$.

Now we can find a reconstruction formula for coefficients of the
series (7).

\begin{teo}
If $u\in H$ then
\begin{equation}
a_n=\lim \limits _{t\to 0}\frac{\Delta ^{(n)}u(t)}{t^{q^n}},\quad
n=0,1,2,\ldots .
\end{equation}
\end{teo}

{\it Proof}. We proceed from the identities
$$
\widetilde{a}^-\left( \lambda \widetilde{f}_n\right) =\lambda
^{1/q}\widetilde{f}_{n-1}\quad n\ge 1,\ \lambda \in \K;\quad
\widetilde{a}^-\widetilde{f}_0=0.
$$
It follows that
$$
\left( \widetilde{a}^-\right) ^nu(t)=\sum \limits _{k=n}^\infty
a_k^{1/q^n}\frac{t^{q^{k-n}}}{D_{k-n}}
$$
whence
$$
a_n^{1/q^n}=\lim \limits _{t\to 0}\frac{\left( \widetilde{a}^-\right)
^nu(t)}t.
$$

Now (8) is a consequence of the identity $\left( \widetilde{a}^-\right)
^n=\sqrt[q^n]{\ }\circ \Delta ^{(n)}$ (see [6]). $\quad \Box $

\section{SMOOTH FUNCTIONS}

Let $u\in C_0(O,\K )$,
$$
\mathfrak D^ku(t)=t^{-q^k}\Delta ^{(k)}u(t),\quad t\in O\setminus
\{ 0\} .
$$
We will say that $u\in C_0^{k+1}(O,\K )$ if $\mathfrak D^ku$ can
be extended to a continuous function on $O$.

$C_0^{k+1}(O,\K )$ can be considered as a Banach space over $\K$,
with the norm
$$
\sup \limits _{t\in O}\left( |u(t)|+|\mathfrak D^ku(t)|\right) .
$$
Note that $C_0^1(O,\K )$ coincides with the set of all
differentiable $\F$-linear  functions $O\to \K$.

In this section we will obtain a characterization of functions
from $C_0^{k+1}(O,\K )$ in terms of coefficients of the expansion
$u=\sum \limits _{n=0}^\infty c_nf_n$. First we prove two
auxiliary results.

\begin{lem}
If a function $v:\ O\setminus \{ 0\} \to \K$ is continuous and
bounded, and $v$ admits a pointwise convergent expansion
\begin{equation}
v(t)=\sum \limits _{n=0}^\infty v_n\tau _n(t),\quad t\in
O\setminus \{ 0\} ,
\end{equation}
$v_n\in \K$, then
$$
\sup \limits _{t\in O\setminus \{ 0\} }|v(t)|=\sup \limits _{0\le
n<\infty }|v_n|.
$$
\end{lem}

{\it Proof}. It is clear that
$$
|v(t)|\le \sup \limits _n|v_n|,\quad t\ne 0.
$$
In order to prove the inverse inequality we will perform the
summation of both sides in (9) through all monic polynomials $t$
with $\deg t=m$, and use the identity
$$
\sum \limits _{{\deg t=m\atop t\mbox{ {\scriptsize monic}}}}g_l(t)G_k(t)=
\left\{
\begin{array}{rl}
0, & \mbox{if }k+l\ne q^m-1,\\
(-1)^m\frac{D_m}{L_m}, & \mbox{if }k+l=q^m-1,
\end{array}\right.
$$
valid for $k<q^m,\ l<q^m$ (see [4, 6]). In particular, we will
need the case when $k=0$, $l=q^n-1$, for which we get
$$
\sum \limits _{{\deg t=m\atop t\mbox{ {\scriptsize monic}}}}g_{q^n-1}(t)=
\left\{
\begin{array}{rl}
0, & \mbox{if }n<m,\\
(-1)^m\frac{D_m}{L_m}, & \mbox{if }n=m.
\end{array}\right.
$$
If $n>m$ then [4, 6] $g_{q^n-1}(t)=t^{-1}e_n(t)=0$ for $\deg t=m$ by the
definition of $e_n(t)$.

Now the summation in (9) yields
$$
\sum \limits _{{\deg t=m\atop t\mbox{ {\scriptsize monic}}}}v(t)=(-
1)^m\frac{D_m}{L_m}\cdot \frac{v_m}{\Gamma _{q^m-1}}.
$$
Using Lemma 1, we obtain that $|v_m|\le
\sup \limits _{t\in O\setminus \{ 0\} }|v(t)|$. $\quad \Box$

\begin{lem}
Let a function $w\in C(O,\K )$ be such that the function $\gamma
(t)=tw(t)$ is $\F$-linear. Then
\begin{equation}
w(t)=\sum \limits _{n=0}^\infty w_n\tau _n(t),\quad w_n\in \K,
\end{equation}
and the series (10) is uniformly convergent on $O$.
\end{lem}

{\it Proof}. Consider an expansion
\begin{equation}
\gamma (t)=\sum \limits \gamma_nf_n(t),\quad t\in O,
\end{equation}
where $\gamma _n\in \K ,\ \gamma _n\to 0$.The fact that the
function $t^{-1}\gamma (t)$ is continuous at $t=0$ means that
$\gamma (t)$ is differentiable at $t=0$. By a result of Wagner
[21], $\gamma _nL_n^{-1}\to 0$ or, equivalently, $|\gamma
_n|q^n\to 0$ for $n\to \infty $.

Dividing both sides of (11) by $t$, we obtain the expansion (10)
with $w_n=D_n^{-1}\gamma _n\Gamma _{q^n-1}$, so that
$w_n\to 0$. $\quad \Box $

Now we are in a position to prove the characterization result.

\begin{teo}
A function $u=\sum \limits _{n=0}^\infty c_nf_n\in C_0(O,\K )$
belongs to $C_0^{k+1}(O,\K )$ if and only if
\begin{equation}
q^{nq^k}|c_n|\to 0\quad \mbox{for }n\to \infty .
\end{equation}
In this case
\begin{equation}
\sup \limits _{t\in O}|\mathfrak D^ku(t)|=\sup \limits
_{n\ge k}q^{(n-k)q^k}|c_n|.
\end{equation}
\end{teo}

{\it Proof}. In view of (6)
$$
\left( a^-\right) ^kf_n=\left\{
\begin{array}{rl}
f_{n-k}, & \mbox{if }n\ge k,\\
0, & \mbox{if }n<k.
\end{array}\right.
$$
Since $\left( a^-\right) ^k=\sqrt[q^k]{\ }\circ \Delta ^{(k)}$, we find
that
\begin{equation}
\mathfrak D^ku(t)=t^{-q^k}\sum \limits _{n=k}^\infty c_nf_{n-k}^{q^k}(t),
\quad t\ne 0,
\end{equation}
which implies the identity
\begin{equation}
\left( \mathfrak D^ku(t)\right) ^{q^{-k}}=\sum \limits _{n=k}^\infty
c_n^{q^{-k}}\Gamma _{q^{n-k}-1}D_{n-k}^{-1}\tau _{n-k}(t),\quad t\ne 0,
\end{equation}

As before, an easy computation shows that $|\Gamma _{q^n-1}D_n^{-1}|=q^n$.
Now, if (12) is satisfied, then by Lemma 2 the right-hand side of
(15) is a continuous function on $O$, which means that
$u\in C_0^{k+1}(O,\K )$. The equality (13) follows from Lemma 3.

Conversely, suppose that $\mathfrak D^ku$ is continuous on $O$.
Let $w(t)=\left( \mathfrak D^ku(t)\right) ^{q^{-k}}$. By (14),
the function $\gamma (t)=tw(t)$ has the form
\begin{equation}
\gamma (t)=\sum \limits _{n=k}^\infty c_n^{q^{-k}}f_{n-k}(t),\quad
t\ne 0.
\end{equation}
Since $w$ is continuous on $O$, we see that $\gamma (0)=0$. On
the other hand, $f_n(0)=0$ for all $n$, so that the equality (16)
holds for all $t\in O$, and the function $\gamma $ is $\F$-linear
(however, we cannot claim the uniform convergence of the series
in (16)).

Now Lemma 4 implies the representation (10) with $w_n\to 0$. It
follows from (10) and (16) that
$$
\sum \limits _{n=0}^\infty
\left(w_n-c_{n+k}^{q^{-k}}L_n^{-1}\right) \tau _n(t)=0
$$
for all $t\ne 0$. By Lemma 3, this means that
$$
w_n=c_{n+k}^{q^{-k}}L_n^{-1}
$$
for $n\ge 0$, which implies (12).  $\Box$

Note that for the case $k=0$ our characterization is equivalent
to Wagner's theorem which was used above in the course of proving
our general result.

\section{ANALYTIC FUNCTIONS}

Denote by $A(O,\K )$ the space of all functions $O\to \K$ of the
form
\begin{equation}
u(t)=\sum \limits _{n=0}^\infty a_nt^{q^n},\quad \K \ni a_n\to 0.
\end{equation}
It is clear that $A(O,\K )\subset C_0(O,\K )$.

\begin{teo}
A function $u=\sum \limits _{n=0}^\infty c_nf_n\in C_0(O,\K )$
belongs to $A(O,\K )$ if and only if
\begin{equation}
q^{\frac{q^n}{q-1}}|c_n|\to 0\quad \mbox{for }n\to \infty .
\end{equation}
\end{teo}

{\it Proof}. Let $u\in A(O,\K )$ be represented by the series
(17). Let us find an expression for its Fourier-Carlitz
coefficients $c_n$.

It is known [6] that
$$
t^{q^n}=\sum \limits _{j=0}^nb_{nj}e_j(t)
$$
where
$$
b_{nj}=\frac{\left. \left( \Delta ^{(j)}t^{q^n}\right) \right|
_{t=1}}{D_j}.
$$
We have seen (Sect. 3) that
$$
\left( \sqrt[q^j]{\ }\circ \Delta ^{(j)}\right) \left(
\frac{t^{q^n}}{D_n}\right) =\left( \widetilde{a}^-\right)
^j\left( \frac{t^{q^n}}{D_n}\right) =\frac{t^{q^{n-j}}}{D_{n-j}}.
$$
whence
$$
\left. \left( \Delta ^{(j)}t^{q^n}\right) \right|
_{t=1}=\frac{D_n}{D_{n-j}^{q^j}}.
$$
so that
$$
t^{q^n}=\sum \limits _{j=0}^n\frac{D_n}{D_{n-j}^{q^j}}f_j(t).
$$

Substituting this into (17) and changing the order of summation
we find that
$$
c_n=\sum \limits _{k=n}^\infty \frac{D_k}{D_{k-n}^{q^n}}a_k.
$$
Since $|D_k|=q^{-\frac{q^k-1}{q-1}}$, we obtain that
$$
|c_n|\le q^{-\frac{q^n-1}{q-1}}\sup \limits _{k\ge n}|a_k|
$$
which implies (18).

Conversely, (18) means that $c_nD_n^{-1}\to 0$ for $n\to \infty
$. Using ($1'$) we can write
\begin{equation}
u(t)=\sum \limits _{n=0}^\infty \frac{c_n}{D_n}\sum \limits
_{j=0}^n(-1)^{n-j}\frac{D_n}{D_jL_{n-j}^{q^j}}t^{q^j}.
\end{equation}
It is easy to get that
$$
\left| \frac{D_n}{D_jL_{n-j}^{q^j}}\right| =q^{s_{nj}}
$$
where
$$
s_{nj}=\left\{
\begin{array}{rl}
(n-j)q^j-q^j-\cdots -q^{n-1}, & \mbox{if }n\ge j+1,\\
0, & \mbox{if }n=j.
\end{array}\right.
$$
Since $s_{nj}\le q^j\left[ (n-j-1)-q^{n-j-1}\right] $ if $n\ge
j+1$, and the function $z\mapsto q^z-z$ increases for $z\ge 1$,
we have $s_{nj}<0$. Now we may rewrite (19) in the form (17) with
$$
|a_n|\le \max \limits _{k\ge n}\left| \frac{c_k}{D_k}\right| \to
0,\quad n\to \infty .\qquad \Box
$$

A similar theorem can be obtained if $K$ is a completion of $\F
(x)$ corresponding to an arbitrary irreducible polynomial $\pi
\in \F [x]$, $\deg \pi =d$. In this case one gets the condition
$$
q^{\frac{q^n}{q^d-1}}|c_n|\to 0\quad \mbox{for }n\to \infty
$$
instead of (18). The proof is based on the fact that the absolute
value $|[i]|_\pi $ equals $q^{-1}$ if $d$ divides $i$, or 1 in
the opposite case (which can be deduced from Lemma 2.13 in [13]).

\section{INDEFINITE SUM}

Viewing the operator $a^-$ as a kind of a derivative, it is
natural to introduce an appropriate antiderivative. Following the
terminology used in the analysis over $\mathbf Z_p$ (see [15]) we
call it the indefinite sum.

Consider in $C_0(O,\K )$ the equation
\begin{equation}
a^-u=f,\quad f\in C_0(O,\K ).
\end{equation}
Suppose that $f=\sum \limits _{k=0}^\infty \varphi _kf_k$,
$\varphi _k\in \K$. Looking for $u=\sum \limits _{k=0}^\infty
c_kf_k$ and using the fact that
$$
a^-u=\sum \limits _{k=1}^\infty c_k^{1/q}f_{k-1}=\sum \limits
_{l=0}^\infty c_{l+1}^{1/q}f_l
$$
we find that $c_{l+1}=\varphi _l^q$, $l=0,1,2,\ldots$. Therefore
$u$ is determined by (20) uniquely up to the term
$c_0f_0(t)=c_0t$, $c_0=u(1)$.

Fixing $u(1)=0$ we obtain a $\F$-linear bounded operator $S$ on
$C_0(O,\K )$, the operator of indefinite sum: $Sf=u$. It follows
from Theorem 2 that $S$ is also a bounded operator on each space
$C_0^k(O,\K )$, $k=1,2,\ldots $. Note that $Sf_k=f_{k+1}$,
$k=0,1,\ldots $, so that $S$ is not compact.

Another possible procedure to find $Sf$ is an interpolation. If
$u=Sf$ then
$$
u(xt)-xu(t)=f^q(t),\quad u(1)=0.
$$
Setting successively $t=1,x,x^2,\ldots$, we get
\begin{align*}
u(x)&=f^q(1), \\
u(x^2)&=f^q(x)+xf^q(1), \\
u(x^3)&=f^q(x^2)+xf^q(x)+x^2f^q(1), \\
......&..............................................\\
u(x^n)&=f^q(x^{n-1})+xf^q(x^{n-2})+\ldots +x^{n-1}f^q(1).
\end{align*}

Since $u$ is assumed $\F$-linear, this determines $u(t)$ for all
$t\in \F [x]$. Extending by continuity, we find $u(t)$ for any
$t\in O$: if
$$
t=\sum _{n=0}^\infty \zeta_nx^n,\quad \zeta_n\in \F \ ,
$$
then
$$
u(t)=\sum _{n=0}^\infty \zeta_nu(x^n)=\sum _{n=1}^\infty \zeta_n
\sum _{j=1}^nx^{j-1}f^q(x^{n-j})
=\sum _{j=1}^\infty \sum _{n=j}^\infty \zeta _nx^{j-1}f^q(x^{n-
j}),
$$
so that
$$
u(t)=\sum _{n=0}^\infty x^nf^q\left( \sum _{m=0}^\infty
\zeta _{m+n+1}x^m\right) .
$$

\section{VOLKENBORN-TYPE INTEGRAL}

The Volkenborn integral of a function on $\mathbf Z_p$ was
introduced in [19] (see also [15]) in order to obtain relations
between some objects of $p$-adic analysis resembling classical
integration formulas of real analysis. Here we extend this
approach to the function field situation. Our definition is based
essentially on the Carlitz difference operator $\Delta $, which
shows again its close connection with basic structures of
analysis over the field $K$.

The integral of a function $f\in C_0^1(O,\K )$ is defined as
$$
\int \limits _{O}f(t)\,dt\stackrel{\mbox{{\footnotesize def}}}{=}\lim \limits
_{n\to \infty }\frac{Sf(x^n)}{x^n}=(Sf)'(0).
$$

It is clear that the integral is a $\F$-linear continuous
functional on $C_0^1(O,\K )$,
$$
\int \limits _{O}cf(t)\,dt=c^q\int \limits _{O}f(t)\,dt,\quad c\in \K .
$$
Since a power series $\sum \limits _{n=0}^\infty a_nt^{q^n}$ with
$a_n\to 0$ converges in $C_0^1(O,\K )$, it can be integrated
termwise:
$$
\int \limits _{O}\sum \limits _{n=0}^\infty a_nt^{q^n}\,dt=
\sum \limits _{n=0}^\infty a_n^q\int \limits _{O}t^{q^n}\,dt.
$$

The integral possesses the following ``invariance'' property
(related, in contrast to the case of $\mathbf Z_p$, to the
multiplicative structure):
\begin{equation}
\int \limits _{O}f(xt)\,dt=x\int \limits _{O}f(t)\,dt-f^q(1).
\end{equation}

Indeed, let $g(t)=f(xt)$. Then
\begin{multline*}
Sg(x^n)=g^q(x^{n-1})+xg^q(x^{n-2})+\cdots +x^{n-1}g^q(1)\\
=f^q(x^n)+xf^q(x^{n-1})+\cdots +x^{n-1}f^q(x)=(Sf)(x^{n+1})-
x^nf^q(1)
\end{multline*}
whence
$$
\frac{Sg(x^n)}{x^n}=x\cdot \frac{(Sf)(x^{n+1})}{x^{n+1}}-f^q(1),
$$
and (21) is obtained by passing to the limit for $n\to \infty $.

Using (20) we obtain by induction that
$$
\int \limits _{O}f(x^nt)\,dt=x^n\int \limits _{O}f(t)\,dt-x^{n-1}f^q(1)
-x^{n-2}f^q(x)-\cdots -f^q(x^{n-1}).
$$
This equality implies the following invariance property. Suppose
that a function $f$ vanishes on all elements $z\in \F [x]$ with
$\deg z<n$. Then, if $g\in \F [x]$, $\deg g\le n$, we have
$$
\int _Of(gt)\, dt=g\int _Of(t)\, dt.
$$

Our last result will contain the calculation of integrals for
some important functions on $O$. Let us recall the definitions of
some special functions introduced by Carlitz (see [8]).

The function $C_s(z),\ s\in \F [x],\ z\in K$, defining the
Carlitz module, is given by the formula
$$
C_s(z)=\sum \limits _{i=0}^{\deg s}f_i(s)z^{q^i}.
$$
If $|z|<1$, $C_s(z)$ can be extended with respect to $s$:
$$
C_s(z)=\sum \limits _{i=0}^\infty f_i(s)z^{q^i},\quad s\in O,
$$
(note that $f_i(s)=0$ for $i>\deg s$ if $s\in \F [x]$). Since
$s^{-1}f_i(s)=L_i^{-1}\tau _i(s)$ and $|L_i^{-1}|=q^i$, the
function $C_s(z)$ is differentiable with respect to $s$.

The Carlitz logarithm $\log _C(z)$ is defined as
$$
\log _C(z)=\sum \limits _{n=0}^\infty (-1)^n\frac{z^{q^n}}{L_n},\quad
|z|<1.
$$

The Carlitz exponential is defined by the power series
$$
e_C(z)=\sum \limits _{n=0}^\infty \frac{z^{q^n}}{D_n},\quad |z|<1.
$$

\begin{teo}
{\rm (i)} For any $n=0,1,2,\ldots$
\begin{equation}
\int \limits _{O}t^{q^n}\,dt=-\frac{1}{[n+1]}.
\end{equation}

{\rm (ii)} For any $n=0,1,2,\ldots$
\begin{equation}
\int \limits _{O}f_n(t)\,dt=\frac{(-1)^{n+1}}{L_{n+1}}.
\end{equation}

{\rm (iii)} If $z\in K,\ |z|<1$, then
\begin{equation}
\int \limits _{O}C_s(z)\,ds=\log _C(z)-z.
\end{equation}
\end{teo}

{\it Proof}. (i) We have seen that
$$
a^-\left( \frac{t^{q^n}}{D_n}\right) =\frac{t^{q^{n-1}}}{D_{n-
1}}, \quad n\ge 1.
$$
It follows from the definition of the operator $S$ that
$$
S\left( \frac{t^{q^n}}{D_n}\right) =\frac{t^{q^{n+1}}-t}{D_{n+
1}}
$$
whence
$$
S\left( t^{q^n}\right) =\frac{D_n^q}{D_{n+1}}\left( t^{q^{n+1}}-
t\right) =\frac{t^{q^{n+1}}-t}{[n+1]},
$$
$$
\frac{S\left( t^{q^n}\right) \left( x^k\right) }{x^k}=
\frac{x^{k(q^{n+1}-1)}-1}{[n+1]}\longrightarrow -
\frac{1}{[n+1]},\quad k\to \infty .
$$

(ii) We have
$$
\int \limits _{O}f_n(t)\,dt=(Sf_n)'(0)=f_{n+1}'(0).
$$
According to ($1'$), the linear term in the expression for $f_i$ is
\begin{equation}
(-1)^i\frac{\left[ i\atop 0\right] }{D_i}t=\frac{(-1)^i}{L_i}t.
\end{equation}
The differentiation yields (23).

(iii) Applying the operator $a^-$ (with respect to the variable
$s$) to the function $C_s(z)$ we find that
$$
a^-_sC_s(z)=\sum \limits _{i=1}^\infty f_{i-1}(s)z^{q^{i-1}}=C_s(z)
$$
whence
$$
S_sC_s(z)=C_s(z)-C_1(z)s=\sum \limits _{i=0}^\infty
f_i(s)z^{q^i}-zs.
$$
Fixing $z$ and denoting $\varphi (s)=S_sC_s(z)$ we obtain that
\begin{equation}
\frac{\varphi (x^n)}{x^n}=\sum \limits _{i=0}^\infty
\frac{f_i(x^n)}{x^n}z^{q^i}-z.
\end{equation}

It is seen from ($1'$) and (25) that
$$
\frac{f_i(x^n)}{x^n}\longrightarrow \frac{(-1)^i}{L_i}\ \
\mbox{for }n\to \infty .
$$

On the other hand,
$$
t^{-1}f_i(t)=D_i^{-1}g_{q^i-1}(t)=D_i^{-1}\Gamma_{q^i-1}\tau
_i(t).
$$
As we know, $|\tau _i(t)|\le 1$, $|D_i^{-1}\Gamma_{q^i-1}|=q^i$,
so that
$$
\left| \frac{f_i(x^n)}{x^n}\right| \le q^i
$$
for all $n$. If $z\in O,\ |z|<1$, then $|z|\le q^{-1}$, and the
series in (26) converges uniformly with respect to $n$. Passing
to the limit $n\to \infty $ in (26) we come to (24). $\quad \Box$

Setting in (24) $z=e_C(t)$, $|t|<1$, we get an identity for the
Carlitz exponential:
$$
\int \limits _Oe_C(st)\,ds=t-e_C(t).
$$
On the other hand, the formula (24) implies a more general
formula (conjectured by D.Goss).

\begin{cor}
If $a\in O,\ z\in K,\ |z|<1$, then
\begin{equation}
\int \limits _OC_{sa}(z)\,ds=a\log _C(z)-C_a(z).
\end{equation}
\end{cor}

{\it Proof}. Since it is shown easily that (for each fixed $z$)
the mapping $O\to C_0^1(O,\K )$ of the form $a\mapsto C_{sa}(z)$
is continuous, it is sufficient to prove (27) for $a=x^n$,
$n=1,2,\ldots $. Using (24), we find that
\begin{equation}
\int \limits _OC_{sx^n}(z)\,ds=x^n(\log _C(z)-z)-\sum \limits
_{k=1}^nx^{n-k}C^q_{x^{k-1}}(z).
\end{equation}

Let $t=\log _C(z)$. Then $z=e_C(t)$ (see [8, 12]). It follows
from properties of $e_C$ [8] that
$$
\sum \limits _{k=1}^nx^{n-k}C^q_{x^{k-1}}(z)=\sum \limits _{k=1}^nx^{n-k}
(e_C(x^kt)-xe_C(x^{k-1}t))=e_C(x^nt)-x^ne_C(t)=C_{x^n}(z)-x^nz.
$$
Substituting this into (28) we come to (27). $\quad \Box$

As $C_{sa}(z)=C_s(C_a(z))$, equation (24) also implies that
$$\int_O C_{sa}(z)\,ds=\log_C(C_a(z))-C_a(z)\,.$$
Comparing this with (27) implies
$$a\log_C(z)=\log_C(C_a(z))$$
which is precisely the functional equation of $\log_C(z)$.

\section*{ACKNOWLEDGEMENTS}

The author is grateful to D.Goss for a lot of valuable remarks
and suggestions, to D.Thakur for a helpful consultation, and to
Z.Yang for pointing out an error in the first version of the
paper. The work was supported in part by the Ukrainian Fund for
Fundamental Research under Grant 1.4/62.

\end{document}